\documentclass[11pt]{article}
\usepackage[T1]{fontenc}
\usepackage{lmodern}
\usepackage[margin=1.05in]{geometry}
\usepackage{amsmath,amssymb,amsthm,mathtools}
\usepackage{microtype}
\usepackage{xcolor}
\usepackage{enumitem}
\usepackage{hyperref}
\hypersetup{
  colorlinks=true,
  linkcolor=blue!55!black,
  citecolor=blue!55!black,
  urlcolor=blue!55!black,
  pdftitle={Quadratic irrational orbits near the middle-third Cantor set: correction and unconditional reductions},
  pdfauthor={Frank Gilson}
}

\newcommand{\Cantor}{\mathcal C}
\newcommand{\dist}{\operatorname{dist}}
\newcommand{\exit}{\operatorname{exit}}
\newcommand{\Q}{\mathbb Q}
\newcommand{\Z}{\mathbb Z}

\numberwithin{equation}{section}
\newtheorem{theorem}{Theorem}[section]
\newtheorem{proposition}[theorem]{Proposition}
\newtheorem{lemma}[theorem]{Lemma}
\newtheorem{corollary}[theorem]{Corollary}
\theoremstyle{definition}
\newtheorem{definition}[theorem]{Definition}
\theoremstyle{remark}

\title{Quadratic Irrational Orbits Near the Middle-Third Cantor Set\\
\large Correction and Unconditional Reductions}
\author{Frank Gilson}
\date{September 10, 2026}

\begin{document}
\maketitle

\begin{abstract}
An earlier form of this manuscript claimed, subject to a ``shallow
contribution'' hypothesis, a log-square upper bound for the first occurrence
of the ternary digit \(1\) in a quadratic irrational, and consequently a
quantitative separation from the middle-third Cantor set. Those conclusions
are withdrawn. The argument used a decomposition defined only after finite
exit had been assumed in order to prove finite exit, and the proposed
quadratic-field-uniform Thue--Mahler estimate did not follow from the cited
\(S\)-unit theory: the relevant set \(S\) varies with the minimal polynomial,
and the displayed \(S\)-unit quotient does not retain the power-of-\(3\)
exponent without an additional normalization.

This replacement records the defects and retains the unconditional part of
the argument. For a quadratic irrational whose orbit exits, we prove the
exact exit-to-distance identity and an elementary height-dependent clearance
bound. We also give a corrected universal bound for a single \(L\)-run and
show that an exceptionally long \(R\)-to-\(L\) transition produces a bounded
value of a binary quadratic norm form at a power of \(3\). A finite-prefix
decomposition identifies the additional uniform orbit-counting theorem that
would be required for a genuine separation result. No assertion that a
quadratic irrational lies outside the Cantor set is made here.
\end{abstract}



\section{Revision notice and current status}

Let \(\Cantor\subseteq[0,1]\) be the middle-third Cantor set.  The earlier
manuscript asserted that, under a stated shallow-orbit hypothesis, a quadratic
irrational \(\alpha\) of naive height \(H\) satisfies
\[
  \exit(\alpha)\ll_{\Q(\alpha)}(\log H)^2
  \quad\text{and}\quad
  \dist(\alpha,\Cantor)\geq H^{-\kappa_{\Q(\alpha)}\log H}.
\]
The former main theorem and both distance corollaries depending on it are
withdrawn and must not be cited as established results.

There are two independent structural defects.

\begin{enumerate}[label=(\arabic*)]
  \item The orbit decomposition used in the shallow-contribution hypothesis
  was defined only under the assumption \(\exit(\alpha)<\infty\).  It was then
  invoked to prove that \(\exit(\alpha)<\infty\).  The same passage contained
  an indexing error: with the convention used here, the first visit to the
  middle interval occurs at time \(\exit(\alpha)-1\), not at time
  \(\exit(\alpha)\).

  \item The claimed quantitative norm-form lemma allowed constants depending
  only on the quadratic field.  In the proposed reduction, however, the set
  of places \(S\) contained the primes dividing \(6a\Delta_f\), where
  \(f(x)=ax^2+bx+c\).  Even when \(\Q(\alpha)\) is fixed, these primes and the
  order discriminant \(\Delta_f\) vary with \(\alpha\).  Standard effective
  bounds for \(S\)-unit equations depend on \(S\), its regulator and prime
  norms, as well as on coefficient heights; see
  \cite{BugeaudGyoryUnit,EvertseGyoryBook,Gyory2019}.  The cited computational
  algorithm of Gherga--Siksek \cite{GhergaSiksek} concerns irreducible binary
  forms of degree at least three and does not directly supply the asserted
  quadratic estimate.
\end{enumerate}

The earlier ancillary archive was also invoked for finite interval checks but
was not included with the manuscript available for this correction.  Those
computer-assisted assertions are therefore not retained below.  The purpose
of this replacement is not to repair the withdrawn theorem by assumption,
but to separate the elementary conclusions that are valid from the uniform
Diophantine and orbit-counting statements that remain unproved.

\section{Dynamical setup}

For \(x\in[0,1)\), put
\[
  \tau(x)=\{3x\},\qquad
  L=[0,1/3),\quad M=[1/3,2/3),\quad R=[2/3,1).
\]
For an irrational \(\alpha\in(0,1)\), write
\[
  \theta_n=\tau^n(\alpha)=\{3^n\alpha\}.
\]
No \(\theta_n\) is a triadic rational, so endpoint conventions do not affect
the orbit.

\begin{definition}[Exit time]
Define
\[
  \exit(\alpha)=
  \inf\{N\geq1:\theta_{N-1}\in M\}
  \in\mathbb N\cup\{\infty\}.
\]
Thus \(\exit(\alpha)\) is the position of the first digit \(1\) in the
ternary expansion of \(\alpha\).  Equivalently,
\(\alpha\in\Cantor\) if and only if \(\exit(\alpha)=\infty\).
\end{definition}

Let \(f(x)=ax^2+bx+c\in\Z[x]\) be the primitive irreducible minimal
polynomial of a quadratic irrational \(\alpha\), and set
\[
  H=H(f)=\max\{|a|,|b|,|c|\}.
\]
Because all rational points used below lie in \([0,1]\), the elementary
estimate
\begin{equation}\label{eq:derivative}
  |f'(x)|=|2ax+b|\leq 3H\qquad(0\leq x\leq1)
\end{equation}
is sufficient.  This improves the unnecessary \(O(H^2)\) derivative bound
used in the earlier version.

\section{Exit time and distance}

The first result is purely dynamical and requires finite exit.

\begin{proposition}[Exact exit-to-distance scaling]\label{prop:exit-distance}
Let \(\alpha\in(0,1)\setminus\Q\) and suppose
\(\exit(\alpha)=N<\infty\).  Then
\[
  \dist(\alpha,\Cantor)
  =3^{-(N-1)}
    \dist\!\left(\theta_{N-1},\{1/3,2/3\}\right).
\]
\end{proposition}

\begin{proof}
The iterates \(\theta_0,\ldots,\theta_{N-2}\) lie in \(L\cup R\), whereas
\(\theta_{N-1}\in M\).  Hence the first \(N-1\) ternary digits of
\(\alpha\) belong to \(\{0,2\}\), and \(\alpha\) lies in a unique
level-\((N-1)\) Cantor interval \(I\).  On \(I\), the map \(\tau^{N-1}\)
is affine with slope \(3^{N-1}\).  The two endpoints of the middle third
removed from \(I\) map to \(1/3\) and \(2/3\).  They are the closest points
of \(\Cantor\) to \(\alpha\), and rescaling their distances gives the
identity.
\end{proof}

For quadratic irrationals, integrality supplies a quantitative clearance at
the exit.

\begin{lemma}[Exit clearance]\label{lem:exit-clearance}
Let \(\alpha\in(0,1)\) be quadratic irrational of height \(H\), and suppose
\(\exit(\alpha)=N<\infty\).  Then
\[
  \dist\!\left(\theta_{N-1},\{1/3,2/3\}\right)
  \geq \frac{1}{9H\,3^N}.
\]
\end{lemma}

\begin{proof}
Let \(m=\lfloor3^{N-1}\alpha\rfloor\) and \(q=3^N\).  Choose
\(p=3m+1\) or \(p=3m+2\), according to which endpoint of \(M\) is nearer
to \(\theta_{N-1}\).  Then \(p/q\in[0,1]\) and
\[
  \left|\alpha-\frac pq\right|
  =3^{-(N-1)}
    \dist\!\left(\theta_{N-1},\{1/3,2/3\}\right).
\]
The integer
\[
  q^2f(p/q)=ap^2+bpq+cq^2
\]
is nonzero, so \(|f(p/q)|\geq q^{-2}\).  By the mean value theorem and
\eqref{eq:derivative},
\[
  q^{-2}\leq |f(p/q)|
  \leq3H\left|\alpha-\frac pq\right|.
\]
Multiplication by \(3^{N-1}=q/3\) yields the claim.
\end{proof}

\begin{corollary}[Unconditional bound in terms of a finite exit]\label{cor:finite-exit}
Under the hypotheses of Lemma~\ref{lem:exit-clearance},
\[
  \dist(\alpha,\Cantor)\geq \frac{1}{3H}\,3^{-2N}.
\]
\end{corollary}

\begin{proof}
Combine Proposition~\ref{prop:exit-distance} and
Lemma~\ref{lem:exit-clearance}:
\[
  \dist(\alpha,\Cantor)
  \geq3^{-(N-1)}\frac{1}{9H3^N}
  =\frac{1}{3H}3^{-2N}.
\]
\end{proof}

Corollary~\ref{cor:finite-exit} does not prove that the exit is finite.  In
particular, it does not address the unresolved possibility that an irrational
algebraic number belongs to \(\Cantor\), which is part of Mahler's problem
\cite{Mahler1984,LevesleySalpVelani}.

\section{A corrected bound for one \texorpdfstring{$L$}{L}-run}

Suppose \(\theta_N\in R\) and write
\[
  \theta_N=\frac23+\delta,qquad0<\delta<\frac13.
\]

\begin{lemma}[Geometry of an \(L\)-run]\label{lem:lrun-geometry}
If \(\theta_{N+1},\ldots,\theta_{N+k}\in L\), then
\[
  \theta_{N+j}=3^j\delta\quad(1\leq j\leq k),
  \qquad 3^k\delta<\frac13.
\]
\end{lemma}

\begin{proof}
The first step is
\(\theta_{N+1}=3(2/3+\delta)-2=3\delta\).  While the orbit remains in
\(L\), the next iterate is obtained by multiplication by \(3\), proving the
formula and its final inequality.
\end{proof}

Define the digit-forced rational
\[
  q_N=3^{N+1},\qquad
  p_N=3\lfloor3^N\alpha\rfloor+2.
\]
Then \(p_N\equiv2\pmod3\), so \(\gcd(p_N,q_N)=1\), and
\begin{equation}\label{eq:forced-approx}
  \alpha-\frac{p_N}{q_N}=\frac{\theta_N-2/3}{3^N}
  =\frac{\delta}{3^N}.
\end{equation}

\begin{lemma}[Integrality bound at an \(R\)-visit]\label{lem:integrality-R}
Let
\[
  A_N=q_N^2f(p_N/q_N)
      =ap_N^2+bp_Nq_N+cq_N^2.
\]
Then \(A_N\in\Z\setminus\{0\}\) and
\[
  1\leq |A_N|\leq 9Hq_N\delta.
\]
Consequently,
\[
  \delta\geq\frac{1}{9H\,3^{N+1}}.
\]
\end{lemma}

\begin{proof}
Integrality is immediate, and nonvanishing follows because \(p_N/q_N\) is
rational while \(f\) is irreducible quadratic.  Both \(\alpha\) and
\(p_N/q_N\) lie in \([0,1]\).  The mean value theorem,
\eqref{eq:derivative}, and \eqref{eq:forced-approx} give
\[
  |A_N|
  \leq q_N^2(3H)\frac{\delta}{3^N}
  =9Hq_N\delta.
\]
Since \(|A_N|\geq1\), the lower bound for \(\delta\) follows.
\end{proof}

\begin{proposition}[Universal single-run bound]\label{prop:single-run}
If \(\theta_N\in R\) and
\(\theta_{N+1},\ldots,\theta_{N+k}\in L\), then
\[
  k<N+2+\log_3H.
\]
\end{proposition}

\begin{proof}
By Lemmas~\ref{lem:lrun-geometry} and \ref{lem:integrality-R},
\[
  \frac{3^k}{9H3^{N+1}}\leq3^k\delta<\frac13.
\]
Thus \(3^k<H3^{N+2}\), which is equivalent to the stated inequality.
\end{proof}

This controls one long excursion through \(L\); it does not bound the number
of such excursions before the orbit enters \(M\).

\section{The valid bounded-norm reduction}

The same calculation yields a useful, but limited, Diophantine reduction.

\begin{proposition}[A deep transition gives a bounded norm value]
\label{prop:bounded-norm}
Suppose \(\theta_N\in R\) and the following \(L\)-run has length \(k\geq1\).
Then
\[
  0<|A_N|\leq9H\,3^{N-k}.
\]
In particular, if for some real \(T\)
\[
  k\geq N+\log_3H+T,
\]
then
\[
  0<|A_N|\leq9\,3^{-T}.
\]
Equivalently, with \(n=N+1\), the coprime pair
\((p_N,3^n)\) satisfies one of the finitely many equations
\[
  ap^2+bp3^n+c3^{2n}=u,
  \qquad 0<|u|\leq\left\lfloor9\,3^{-T}\right\rfloor.
\]
\end{proposition}

\begin{proof}
Lemma~\ref{lem:lrun-geometry} gives \(\delta<3^{-(k+1)}\).  Substitution in
Lemma~\ref{lem:integrality-R} gives
\[
  |A_N|\leq9H3^{N+1}3^{-(k+1)}=9H3^{N-k}.
\]
The remaining statements follow immediately.
\end{proof}

Completing the square rewrites the displayed equation as
\[
  (2ap+b3^n)^2-\Delta_f3^{2n}=4au,
  \qquad \Delta_f=b^2-4ac.
\]
For a fixed polynomial \(f\) and fixed \(u\), this is a Pell/norm-form type
equation with a power-of-\(3\) restriction.  Proposition~\ref{prop:bounded-norm}
does not provide a bound for \(n\), uniformly or otherwise.

The earlier manuscript attempted to obtain such a bound by factoring in
\(K=\Q(\sqrt{\Delta_f})\).  If
\[
  z=2ap+(b+\sqrt{\Delta_f})3^n=\gamma\varepsilon,
\]
one obtains an \(S\)-unit equation using
\[
  x=\frac{\varepsilon'}{\varepsilon},
  \qquad y=\frac{3^n}{\varepsilon}.
\]
However, the simultaneous replacement
\[
  (n,\varepsilon)\longmapsto(n+m,3^m\varepsilon)
\]
leaves both \(x\) and \(y\) unchanged.  A height bound for the resulting
\(S\)-unit solution therefore does not by itself bound \(n\).  The coprimality
and integrality conditions may contain additional information, but a valid
argument must preserve and use them.  No such uniform normalization is
proved here.

\section{A noncircular finite-prefix decomposition}

The original decomposition can be stated correctly without assuming a
finite exit.  This formulation makes the missing orbit-counting input
transparent.

\begin{lemma}[Decomposition of an \(M\)-avoiding prefix]
\label{lem:prefix-decomposition}
Let \(J\geq N_*\geq0\).  Suppose
\[
  \theta_j\notin M\quad(0\leq j\leq J),
  \qquad \theta_{N_*}\in R.
\]
Let
\[
  \mathcal R_J=\{n:N_*\leq n\leq J\text{ and }\theta_n\in R\}.
\]
For \(n\in\mathcal R_J\), define the truncated subsequent \(L\)-run length
\[
  \ell_J(n)=\max\{\ell:0\leq\ell\leq J-n,
      \ \theta_{n+1},\ldots,\theta_{n+\ell}\in L\}.
\]
Then the integer intervals
\[
  I_n=\{n,n+1,\ldots,n+\ell_J(n)\}
  \qquad(n\in\mathcal R_J)
\]
are pairwise disjoint and partition \(\{N_*,\ldots,J\}\).  Hence
\[
  J-N_*+1
  =\sum_{n\in\mathcal R_J}\bigl(1+\ell_J(n)\bigr).
\]
\end{lemma}

\begin{proof}
Every \(R\)-time \(n\) begins its interval \(I_n\), followed by the maximal
run of \(L\)-times within the prefix.  The next \(R\)-time, if present, occurs
strictly after that interval, so the intervals are disjoint.  Conversely,
given \(t\in[N_*,J]\), either \(\theta_t\in R\), in which case \(t\in I_t\),
or \(\theta_t\in L\).  Moving backward through the current \(L\)-run reaches
an \(R\)-time no earlier than \(N_*\), because \(M\) is absent and
\(\theta_{N_*}\in R\).  Thus \(t\) lies in the interval beginning at that
\(R\)-time.
\end{proof}

Lemma~\ref{lem:prefix-decomposition} shows exactly what a future separation
proof must add.  To force a finite exit, one needs a bound, independent of
the cutoff \(J\), for the right-hand side whenever the orbit avoids \(M\)
through time \(J\).  Proposition~\ref{prop:single-run} bounds the length of
each individual run in terms of its starting time, not the number or total
contribution of the runs.  Proposition~\ref{prop:bounded-norm} only reduces a
restricted class of unusually deep transitions to norm-form equations.
Neither result supplies the required uniform total bound.

\begin{proposition}[A correct conditional implication]
\label{prop:conditional-template}
Fix a class \(\mathcal A\) of quadratic irrationals.  Suppose there are
nonnegative integer-valued functions \(B(H)\) and \(E(H)\) such that, for every
\(\alpha\in\mathcal A\) of height \(H\):
\begin{enumerate}[label=(\roman*)]
  \item the orbit reaches \(M\cup R\) by time \(E(H)\); and
  \item whenever \(N_*\leq J\), \(\theta_{N_*}\in R\), and the orbit avoids
  \(M\) through time \(J\), one has
  \[
    \sum_{n\in\mathcal R_J}(1+\ell_J(n))\leq B(H).
  \]
\end{enumerate}
Then every \(\alpha\in\mathcal A\) has finite exit and
\[
  \exit(\alpha)\leq E(H)+B(H)+1.
\]
\end{proposition}

\begin{proof}
If the first visit to \(M\cup R\) is a visit to \(M\), the conclusion is
immediate.  Otherwise let \(N_*\leq E(H)\) be the first \(R\)-time.  If the
orbit avoided \(M\) through time \(J=N_*+B(H)\),
Lemma~\ref{lem:prefix-decomposition} would give
\(B(H)+1\leq B(H)\), a contradiction.  Hence \(M\) is reached by time
\(N_*+B(H)\), and the digit-index definition adds one to obtain the stated
exit bound.
\end{proof}

Taking \(B(H)=O((\log H)^2)\) in
Proposition~\ref{prop:conditional-template} would recover the shape of the
withdrawn exit estimate.  But condition (ii) is precisely the missing
uniform orbit-counting statement; labelling it a hypothesis does not
establish a separation theorem.

\section{Concluding status}

The unconditional conclusions of this replacement are:
\begin{enumerate}[label=(\arabic*)]
  \item the exact distance formula of Proposition~\ref{prop:exit-distance}
  and the finite-exit lower bound of Corollary~\ref{cor:finite-exit};
  \item the corrected bound \(k<N+2+\log_3H\) for one \(L\)-run; and
  \item the bounded quadratic norm-form identity of
  Proposition~\ref{prop:bounded-norm} for a sufficiently deep transition.
\end{enumerate}

These facts do not imply that a quadratic irrational exits the Cantor
itinerary.  A genuine quantitative separation result would require both a
noncircular uniform estimate for all \(M\)-avoiding prefixes and a correctly
parameterized treatment of the associated quadratic norm-form equations.
The dependence on the varying set of bad primes and the rational
power-of-\(3\) scaling must be controlled explicitly.  Those problems are
left open.

\end{document}